\documentclass[a4paper]{article}

\usepackage{graphicx}          
 \usepackage{rotating}
\usepackage{amsmath,amssymb}
\usepackage{diagbox}
 \usepackage{lipsum}
\usepackage{epstopdf}
\usepackage{color}

\usepackage{booktabs}
 \usepackage{multirow}
\usepackage{subcaption}
\usepackage[titletoc]{appendix}

\usepackage[tikz]{bclogo}
\usepackage{multicol}
\usepackage{mdframed}
\mdfdefinestyle{MyFrame}{%
    linecolor=black,
    outerlinewidth=1pt,
    roundcorner=5pt,
    innertopmargin=\baselineskip,
    innerbottommargin=\baselineskip,
    innerrightmargin=10pt,
    innerleftmargin=10pt,
    backgroundcolor=green!10!white}
\usepackage{cite}

\usepackage{scalerel}
\usepackage[normalem]{ulem}

\newtheorem{corollary}{Corollary}
\newtheorem{lemma}{Lemma}
\newtheorem{problem}{Problem}

\newtheorem{property}{Property}

\def\c{\boldsymbol{c}}

\title{On Dubins Paths to a Circle}

\author{Zheng Chen\thanks{Professor, School of Aeronautics and Astronautics, Zhejiang University. Email: z-chen\@ zju.edu.cn}\\
Zhejiang University, Hangzhou 310027, Zhejiang, China}

\begin{document}

\maketitle{}

\begin{abstract}
This paper is concerned with characterizing  the shortest path of a Dubins vehicle  from a position with a prescribed  heading angle  to a target circle with the final heading tangential to the target circle. Such a shortest path is of significant  importance as it is usually required in real-world scenarios, such as taking a snapshot of a forbidden region or  loitering above a ground sensor to collect data by a fixed-wing unmanned aerial vehicle in a minimum time.
By applying Pontryagin's maximum principle, some geometric properties for the shortest path are  established {without considering any assumption on the relationship between the minimum turning radius and radius of the target circle}, showing that the shortest path must lie in a sufficient family of  12 types. {
By employing those geometric properties, the analytical solution to each type is devised so that the length of each type can be computed in a constant time. In addition, some properties depending on problem parameters are found so that 
the shortest path can be computed without checking all the 12 types.}  Finally, some numerical simulations are presented, illustrating and validating the developments of the paper.
\end{abstract}



\section{Introduction}

When  planning a minimum-time path for an Unmanned Aerial Vehicle (UAV), it is common to consider that the UAV flies in altitude hold mode and that its cruise speed  is constant. If taking into account the  constraint that the turning rate is bounded from below, the UAV can be considered as a typical nonholonomic vehicle that  moves only forward at a constant speed with a minimum turning radius. Such a nonholonomic vehicle has been {named}  Dubins vehicle since Dubins studied its shortest path in 1957 \cite{Dubins:57}. For this reason, the problem of minimum-time path planning for a large class of vehicles, such as UAVs \cite{Matveev:11}, fixed-wing aircrafts \cite{Lugo:14}, and thrusted skates \cite{Lynch:03},  is simplified in the literature to the problem of finding the shortest path of Dubins vehicles (note that the shortest path is equivalent to the minimum-time path as the speed is constant).

By geometric analysis, it is shown in \cite{Dubins:57} that the shortest path of Dubins vehicle between two configurations (a configuration consists of a position point and a heading orientation angle) is a $C^1$ path which is a concatenation of circular arcs and straight line segments, and that it can be computed in a constant time by checking at most 6 candidates. This result was proved later in \cite{Boissonnat:94,Sussmann:94} using the optimal control theory. In \cite{Bui:94}, the shortest Dubins path from a configuration to a point with a free terminal heading angle was studied as well, and this problem is now {named} the relaxed Dubins problem. In addition, the reachability domains of Dubins vehicle in some scenarios have been studied recently \cite{Berdyshev:06,Fedotov:18,Patsko:18,Patsko:19}. 

Following the aforementioned papers,  the shortest Dubins paths in many scenarios  have been studied in the literature.  For instance, the shortest Dubins paths through three points were studied  in \cite{Chen:19Automatica,Ma:06,Sadeghi:16,Chen:med19}. To be specific, Refs.~\cite{Ma:06,Chen:med19} presented a natural extension of the relaxed Dubins problem in \cite{Bui:94}, which consists of moving from a configuration, via an intermediate point, to a fixed  final point with free heading angles at both the intermediate and final points. In \cite{Sadeghi:16}, the three-point Dubins problem (consisting of three points with prescribed heading angles at initial and final points) was studied with an assumption that the distance between any two consecutive points was {greater than} four times  the minimum turning radius. By removing the assumption in \cite{Sadeghi:16}, the three-point Dubins problem was thoroughly studied in \cite{Chen:19Automatica}, showing that the shortest path of three-point Dubins problem must be in a sufficient family of 18 types, and a polynomial-based method was proposed to efficiently compute each of the 18 types.
Based on the solutions established in \cite{Chen:19Automatica,Ma:06,Sadeghi:16}, some algorithms have been developed to address the Dubins traveling salesman problem  \cite{Isaiah:15}.



In all the aforementioned papers, the shortest Dubins paths were studied in scenarios where the Dubins vehicle  moves between or among fixed points. In this paper, {another significant scenario that the Dubins vehicle moves from a  configuration to a target circle with the terminal velocity tangential to the circle, will be studied.} This study is motivated by surveillance missions requiring a UAV to take a snapshot of an adversary radar or to inspect a forbidden region.  There are some other applications that require planning a minimum-time path from a configuration to a target circle with the terminal heading tangential to the circle. For instance as proposed in \cite{Manyam:19},  if a UAV collects data from a sensor on the ground, it has to loiter above a sensor while a reliable communication network is established. To loiter above a ground sensor, the UAV should travel along a circle with smallest radius possible with the center of the circle above the sensor location. In addition, if a UAV is to patrol a circular-like border from an airport, the same scenario appears because the UAV needs to fly from an initial configuration to a circle with the heading tangential to the circle \cite{Matveev:11}. As pointed out in \cite{Manyam:19}, another application that is closely related to this problem is the Dubins traveling salesman problem with neighborhoods  \cite{GuimaraesMacharet:2012,Vana:15,Isaacs:11,Isaacs:13}.



{In fact, such a problem has been studied in a seminal work \cite{Balluchi:1996}, which presents some geometric properties for the solution of the problem with an assumption that the minimum turning radius is the same as that of target circle. Recently, with a strict assumption that the distance between the initial position and the target circle is greater than four times the minimum turning radius, the shortest Dubins path to a circle has been  studied in \cite{Manyam:19}.} In this paper, the optimal control theory will be used to characterize the shortest Dubins path to a circle  without any assumptions. First of all, the shortest Dubins path to a circle is formulated as the solution path of an optimal control problem. Then, by using Pontryagin's Maximum Principle (PMP) \cite{Pontryagin}, the necessary conditions for optimality are synthesized so that some geometric properties for the solution path are established, showing that the solution path is a smooth concatenation of circular arcs and straight line segments. Further analyses of those geometric properties not only show that some results in \cite{Manyam:19} are not correct but also 
restrict the solution path into a sufficient family of 12 types. These geometric properties are also used so that the analytical solution to each of all the 12 types is devised. {In addition, some geometric properties depending on problem parameters are found so that some of the 12 types can be ruled out for some specific initial conditions.}   As a result, the shortest Dubins path from any initial configuration to any target circle  can be computed in a constant time by checking the values of at most 12 analytical functions.





\section{Problem formulation}\label{SE:Problem}

In this section, the optimal control problem is formulated and its necessary conditions are derived from the PMP. 

\subsection{Optimal control problem}

\begin{figure}[!htp]
\centering
\includegraphics[width = 3cm]{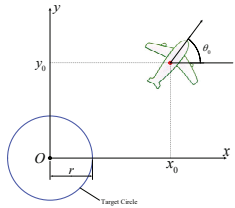}
\caption{Geometry for the engagement.}
\label{Fig:geometry}
\end{figure}

Without lose of generality, consider that the coordinate system $Oxy$ has its origin located at the center of the target circle with radius $r> 0$, as presented in Fig.~\ref{Fig:geometry}. The state of the UAV in frame $Oxy$ consists of a position vector and a heading orientation angle. We denote by $(x,y)\in \mathbb{R}^2$ the position of the UAV, and by $\theta\in [0,2\pi]$ the heading orientation angle of the UAV with respect to the $x$ axis, positive when measured counter-clockwise.

The kinematic  model assumes that the speed is constant and that the turning radius is lower bounded. By normalizing the position $(x,y)$ so that the speed is one, the equations of motion  can be expressed as
\begin{align}
(\Sigma):\ \ \ \ 
\begin{cases}
\dot{x}(t) \ = \cos \theta(t),\\
\dot{y}(t) \ = \sin \theta(t),\\
\dot{\theta}(t)\  = u(t)/\rho
\end{cases}
\label{EQ:Sigma}
\end{align}
where $t\geq 0$ is time,  $u\in[-1,1]$ is the control input, and  $\rho>0$ denotes the minimum turning radius of the UAV.  The initial state is  $ (x_0,y_0,\theta_0)$.

Define two real-valued functions as
\begin{align}
 \phi_1(x,y,\theta) \ &= \frac{1}{2}(x^2 + y^2 - r^2)\\
\phi_2(x,y,\theta) \ &= x\cos(\theta) + y \sin(\theta)
\end{align}
Then, the minimum-time path of the UAV  from  $(x_0,y_0,\theta_0)$ to the target circle with the terminal heading tangential to the circle is equivalent to the solution path of the following Optimal Control Problem (OCP).
\begin{problem}[OCP]
The OCP consists of finding a minimum time $t_f > 0$  so that there exists a path $(x(\cdot),y(\cdot),\theta(\cdot))$ over $[0,t_f]$ subject to $(\Sigma)$ and $u\in[-1,1]$ with boundary conditions:
$$(x(0),y(0),\theta(0)) = (x_0,y_0,\theta_0)$$
and 
\begin{align}
\phi_1(x(t_f),y(t_f),\theta(t_f)) = 0\\
\phi_2(x(t_f),y(t_f),\theta(t_f)) = 0
\end{align}
\end{problem}

If a UAV is required to finally move along the target circle, we have $r\geq \rho$.   {For completeness}, the cases of $\rho > r$, $\rho = r$, and $\rho < r$ are all considered in the paper. 


\subsection{Necessary conditions}
Denote by $p_x$, $p_y$, and $p_{\theta}$ the costates of $x$, $y$, $\theta$, respectively. {Since the goal is to maximize $\int_0^{t_f} -1 dt$, it follows from \cite{Chen:19Automatica} that} the Hamiltonian  is 
\begin{align}
H = p_x \cos \theta + p_y \sin \theta + p_{\theta} u/\rho -1.
\label{EQ:Hamiltonian}
\end{align}
According to PMP \cite{Pontryagin}, for $t\in [0,t_f]$ it holds that 
\begin{align}
\dot{p}_x(t)\ & = -\frac{\partial H}{\partial x} = 0,\label{EQ:dp-x}\\
\dot{p}_y(t)\ &  = - \frac{\partial H}{\partial y} = 0,\label{EQ:dp-y}\\
\dot{p}_{\theta}(t)\ & = -\frac{\partial H}{\partial \theta} = p_x(t) \sin \theta(t) - p_y(t) \cos \theta(t).\label{EQ:dp-theta}
\end{align}
As Eq.~(\ref{EQ:dp-x}) and Eq.~(\ref{EQ:dp-y}) indicate that $p_x$ and $p_y$ are constant, we hereafter use $p_x$ and $p_y$ to denote $p_x(t)$ and $p_y(t)$, respectively.  Integrating Eq.~(\ref{EQ:dp-theta}) leads to
\begin{align}
p_{\theta}(t) = p_x y(t) - p_y x(t) + c_0,
\label{EQ:p-theta}
\end{align}
where $c_0\in\mathbb{R}$ is a constant. If $p_{\theta}\equiv 0$ on a nonzero interval,   Eq.~(\ref{EQ:p-theta}) implies that  the {path} $(x,y)$ on this interval forms a straight line, further indicating $u \equiv 0$ on this interval. Thus, the PMP implies
\begin{align}
u=
\begin{cases}
1,\ \ \ \ p_{\theta} > 0,\\
0,\ \ \ \ p_{\theta}\equiv 0,\\
-1,\ \ p_{\theta}<0.
\end{cases}
\label{EQ:control}
\end{align}
 The transversality conditions are
\begin{eqnarray}
p_x  =\ \frac{\partial \phi_1}{\partial x} \nu_1 + \frac{\partial \phi_2 }{\partial x}\nu_2  = \nu_1 x(t_f) + \nu_2 \cos\theta(t_f)\label{EQ:trans1}\\
p_y   =\ \frac{\partial \phi_1}{\partial y} \nu_1 + \frac{\partial \phi_2 }{\partial y}\nu_2 = \nu_1 y(t_f) + \nu_2 \sin\theta(t_f)\label{EQ:trans2}\\
p_{\theta} (t_f)  =\ \frac{\partial \phi_1}{\partial \theta} \nu_1 + \frac{\partial \phi_2 }{\partial \theta}\nu_2 \ \ \ \ \ \ \ \ \ \ \  \ \ \ \ \ \ \ \ \ \ \ \ \ \ \ \ \ \ \ \ \   \nonumber\\
 =  \nu_2\{ -x(t_f)\sin[\theta(t_f)] + y(t_f) \cos[\theta(t_f)]\}\label{EQ:trans3}
\end{eqnarray}
where $\nu_1$ and $\nu_2$ are Lagrangian multipliers. 


\section{Characterization of the shortest path}\label{SE:Property}


 %


 Denote by ``S'' and ``C'' a straight line segment and a circular arc with   radius of $\rho$, respectively. We immediately have the following property.
 

\begin{property}\label{LE:path-type}
The solution path of the OCP is of type $CCC$ or $CSC$ or a substring thereof, where
\begin{itemize}
\item CCC = \{RLR,\ LRL\},
\item CSC = \{RSR, \ RSL,\ LSR,\ LSL\}. 
\end{itemize}
where ``R'' and ``L''  means the corresponding circular arcs have right and left turning directions. respectively.
\end{property} 
\noindent {\it Proof}. Let $\gamma(\cdot)$ on $[0,t_f]$ be  the solution path of the OCP.   By contradiction, assume that $\gamma(t)$ does not belong to either CCC or CSC or a substring thereof. Then, according to \cite{Dubins:57,Sussmann:94}, there is a shorter path from $\gamma(0)$ to $\gamma(t_f)$,  contradicting with the assumption that $\gamma(t)$ is the shortest. $\square$


\begin{lemma}[Balluchi and Sou\'eres \cite{Balluchi:1996}]\label{LE:S-pass-center}
If solution path of the OCP is of CSC, then {the straight line segment S is collinear with the center of the target circle.}
\end{lemma}
\noindent {\it Proof}.
 Eliminating  $\nu_1$ and $\nu_2$ in Eqs.~(\ref{EQ:trans1}--\ref{EQ:trans3}) leads to
\begin{align}
p_{\theta}(t_f) = p_x y(t_f) - p_y x(t_f).
\label{EQ:p-theta-tf}
\end{align}
Combining  Eq.~(\ref{EQ:p-theta-tf}) with Eq.~(\ref{EQ:p-theta}) indicates $c_0 = 0$ and
\begin{align}
p_{\theta}(t) = p_x y(t) - p_y x(t),\ \forall t\in[0,t_f].
\label{EQ:p-theta-t}
\end{align}
Assume that $[t_1,t_2]\subset [0,t_f]$ is the interval of the straight line segment S. Then, since $p_{\theta}\equiv 0$ along  S, it follows
\begin{align}
 p_x y(t) - p_y x(t) = 0,\ \forall t\in[t_1,t_2].
\nonumber
\end{align}
Note that the the center of target circle (or the origin of $Oxy$) is colinear with  the straight line $p_x y - p_y x = 0$,  completing the proof. $\square$ 

\begin{corollary}\label{CO:1}
The solution path of the OCP must be one of types in 
$
\mathcal{F} = \{CCC,CSC,SC,CC,C\}
$,
where 
\begin{align}
 CCC\ & = \{ RLR, LRL \},\nonumber\\
 CSC \ &= \{RSR, RSL, LSL, LSR\},\nonumber\\
SC \ &= \{SR,SL\},\nonumber\\
CC \ &= \{RL,LR\},\nonumber\\
C \ &= \{R, L\}.\nonumber
\end{align}
\end{corollary}
\noindent {\it Proof.} Lemma \ref{LE:S-pass-center} indicates that the solution of OCP can never be of type CS, completing the proof. $\square$

\noindent Note that the total number of possible types (including substrings) for the shortest Dubins path between two configurations is up to 15 \cite{Dubins:57,Sussmann:94}. However, this corollary means that, including the substrings, the total number of possible types for solution paths of the OCP is 12.  

\begin{figure}[htp]
\centering
\begin{subfigure}[t]{4cm}
\centering
\includegraphics[width =4cm]{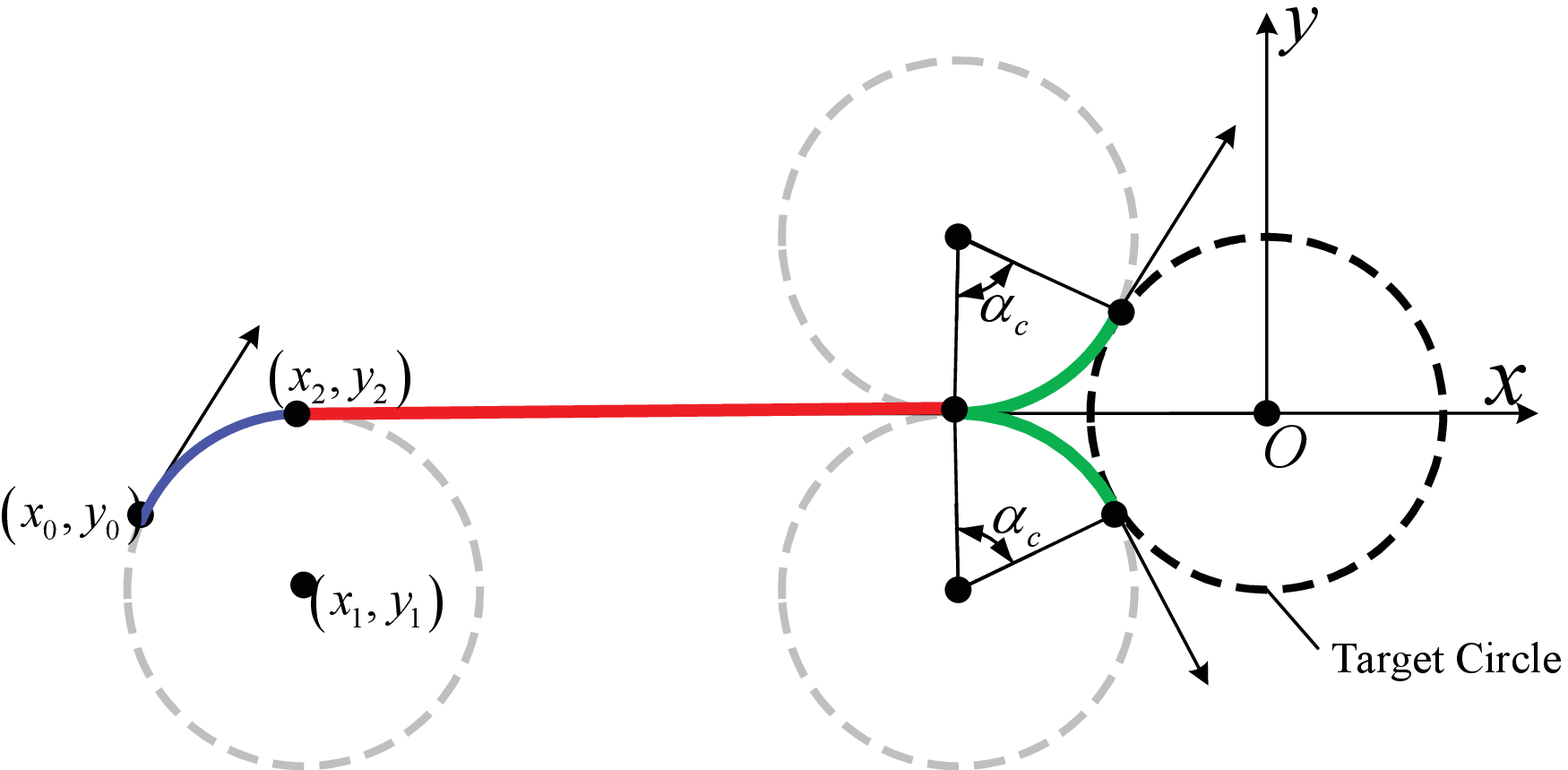}
\caption{RSC}
\label{Fig:RSC}
\end{subfigure}
~~~
\begin{subfigure}[t]{4cm}
\centering
\includegraphics[width =4cm]{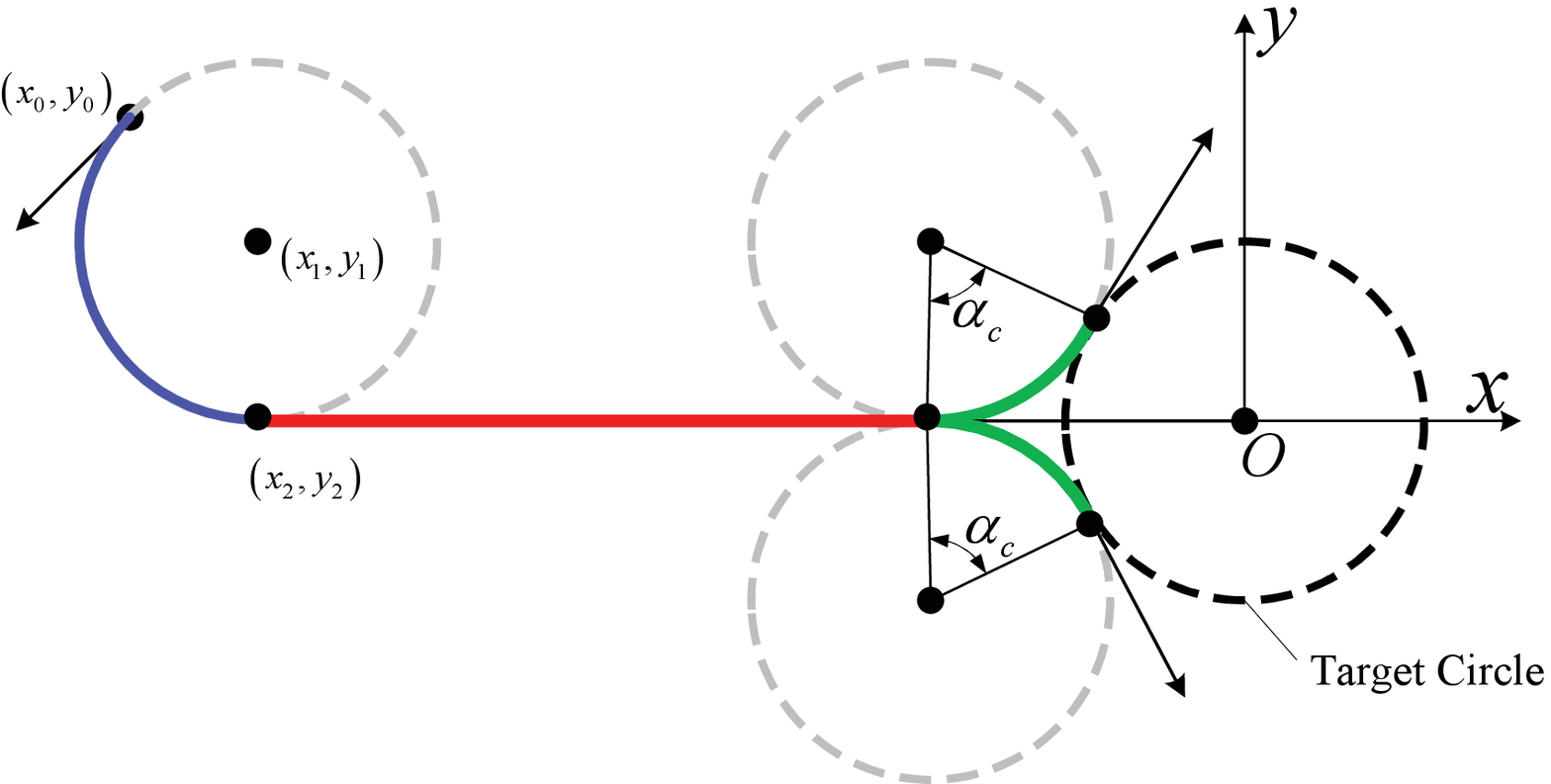}
\caption{LSC}
\label{Fig:LSC}
\end{subfigure}
\caption{The geometry for CSC with the final circular arc externally tangent to the target circle.}
\label{Fig:CSC}
\end{figure}

\begin{lemma}\label{LE:alpha}
Assume that the solution path of the OCP is of type $C_1S_2C_3$, and let $\alpha_c>0$ be the radian of the final circular arc $C_3$. Then, the following two statements hold:
\begin{description}
\item (1) If $C_3$ is externally tangent to the target circle, we have
\begin{align}
\alpha_c = \arccos({\rho}/{\rho + r}).
\label{EQ:lemma-alphac}
\end{align}
\item (2) If $C_3$ is inside the target circle, we have
\begin{align}
\alpha_c = \pi - \arccos({\rho}/{r- \rho}).
\label{EQ:lemma-alphac1}
\end{align}
\end{description}
\end{lemma}
\noindent Proof. (1) According to Lemma \ref{LE:S-pass-center},  $C_3$ must be tangent to both  $S$ and the target circle, as illustrated in Fig.~\ref{Fig:CSC}. 
According to the geometry, the radian $\alpha_c>0$ of  $C_3$ takes a value such that 
$\cos\alpha_c = \frac{\rho}{r+\rho}.$
Since $\alpha_c < \pi$, Eq.~(\ref{EQ:lemma-alphac}) holds, completing the proof of the first statement.

\begin{figure}
\centering
\includegraphics[width = 1.5 in]{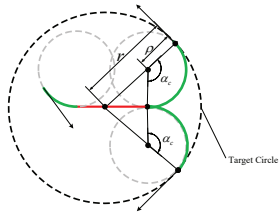}
\caption{The geometry for LSC paths with the final circular arc internally tangent to the target circle.}
\label{Fig:Internal_tangent}
\end{figure}

(2) Analogously, $C_3$ must be tangent to both $S$ and the target circle, as illustrated  in Fig.~\ref{Fig:Internal_tangent}.
The geometry in Fig.~\ref{Fig:Internal_tangent} implies that the radian $\alpha_c>0$ of  $C_3$ takes a value such that 
$\cos (\pi - \alpha_c) = \frac{\rho}{r - \rho}.$
Since $\alpha_c < \pi$, we  have Eq.~(\ref{EQ:lemma-alphac1}), completing the proof. $\square$


 \begin{lemma}[Balluchi and Sou\'eres \cite{Balluchi:1996}]\label{LE:CCC-Straight}
If the solution path of the OCP is of type $C_1C_2C_3$, then the concatenating point from $C_1$ to $C_2$,  the concatenating point from $C_2$ to $C_3$,  and the center of the target circle  are collinear.
\end{lemma}
\noindent Proof. Let $(x_1,y_1)\in\mathbb{R}^2$ be the concatenating point from $C_1$ to $C_2$, and let $(x_2,y_2)\in\mathbb{R}^2$ be the concatenating point from $C_2$ to $C_3$. 
According to Eq.~(\ref{EQ:control}),  we have $p_{\theta}>0$ (resp. $p_{\theta} < 0$) on any left (resp. right) turning circular arc.  Note that $CCC$ can be either RLR or LRL. Because the costate variable $p_{\theta}$ is continuous, we have that $p_{\theta} = 0$ at  $(x_1,y_1)$ and $(x_2,y_2)$. As  $c_0 = 0$ from Eq.~(\ref{EQ:p-theta-t}), writing $p_{\theta}$ in Eq.~(\ref{EQ:p-theta}) explicitly at the two points  $(x_1,y_1)$ and $(x_2,y_2)$ leads to
 $$ p_x y_1 - p_y x_1 = p_x y_2 - p_y x_2 =  0.$$
Hence, the two  points $(x_1,y_1)$ and $(x_2,y_2)$ lie on the straight line of $p_x y - p_y x = 0$. Note that the center of the target circle (the origin of frame $Oxy$) also lies on the straight line of $p_x y - p_y x = 0$, completing the proof. $\square$
 


\begin{lemma}\label{LE:3}
If the solution path is of type $C_1S_2C_3$, the following two statements hold:
\begin{description}
\item (1) If $\rho \geq r/2$,  the circular arc $C_3$ and the target circle are externally tangent to each other.
\item (2) If $\rho < r/2$, the center of the target circle does not lie in the interior of the circle coinciding with $C_1$. 
\end{description}
\end{lemma}
\noindent Proof. (1) By contradiction, assume that   $C_3$ and the target circle are internally tangent to each other. If $\rho>r$, the contradicting assumption means that the target circle lies inside $C_3$. Then, the velocity along $S_2$ cannot point to the center of the target circle, controdicting with Lemma \ref{LE:S-pass-center}. Hence, if $\rho>r$, this lemma holds.  If $\rho = r$, the contradicting assumption implies that $C_3$ coincides with the target circle. In this case, the velocity along $S_2$ will not point to center of the target circle, contradicting with Lemma \ref{LE:S-pass-center}. Hence, if $r = \rho$, this lemma holds. If  $r/2 \leq \rho < r $, the contradicting assumption indicates that the center of the target circle lies inside $C_3$. Analogously, the velocity along $S_2$ cannot point to center of the target circle,  contradicting with Lemma \ref{LE:S-pass-center}. Hence,  this lemma holds if $r/2 \leq \rho < r $.  

(2) By contradiction, assume that the center of target circle lies in the interior of $C_1$.  Then, $S_2$ cannot be collinear with the center of the target circle. This contradicts with Lemma \ref{LE:S-pass-center}, completing proof. $\square$


\begin{lemma}\label{LE:CCC-ruleout-internal}
Assume $\rho > r$. A path of type $C_1C_2C_3$ cannot be optimal if  the target circle is internally tangent to $C_3$.
\end{lemma}

\begin{figure}[!htp]
\centering
\includegraphics[width = 1.5in]{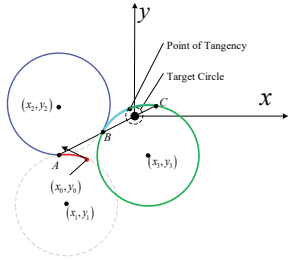}
\caption{An LRL path with the target circle lying inside the final circular arc.}
\label{Fig:LRL_internal}
\end{figure}

\noindent Proof.  If the target circle is internally tangent to $C_2$, the assumption of $\rho > r$ indicates the target circle lies inside $C_3$, as presented by the LRL path in Fig.~\ref{Fig:LRL_internal}.
By contradiction, assume that $C_1C_2C_3$ is optimal. According to Lemma \ref{LE:CCC-Straight}, the point $A$, the point $B$, and the origin  are collinear, implying that  the final point must be on the circular arc from $C$ to $B$ counter-clockwise. According to \cite[Lemma 3]{Bui:94}, the path in this case is not the shortest. Hence, by contraposition, this lemma holds if the path is LRL. As for  RLR, we can prove this lemma by the same way.
$\square$

\begin{lemma}
A path of type $C_1C_2C_3$ cannot be optimal if the distance from the center of $C_1$ to the center of the target circle is not less than $r+\rho$.
\end{lemma}
\noindent Proof. Denote by $d>0$ the distance from the center of $C_1$ to the center of the target circle. We first consider $d = r+\rho$. Then, the solution path is a single circular arc instead of a type $CCC$. Hence, by contraposition, a path of type $C_1C_2C_3$ is not optimal if $d=r+\rho$. Next, we assume $d>r+\rho$. In this case, there exists a straight line segment that is tangent to both $C_1$ and the target circle so that a path of $CS$ is shorter than $C_1C_2C_3$. Hence, a path of type $C_1C_2C_3$ is not optimal if $d > r+\rho$, completing the proof.
$\square$

\begin{figure}[!htp]
\centering
\includegraphics[width = 4cm]{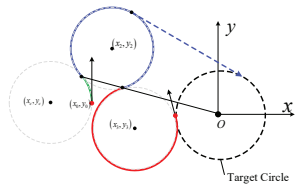}
\caption{An LRL path with the target circle externally tangent to the final circular arc.}
\label{Fig:CCC_external_tangent}
\end{figure}

\begin{lemma}\label{LE:7}
Given a path of type $C_1C_2C_3$, denote by $d_1>0$ and $d_2>0$ the distances from the center of target circle to the centers of $C_1$ and $C_2$, respectively. If $d_1> \rho + r$ and $d_2 > \rho + r$, the path of $C_1C_2C_3$ cannot be optimal. 
\end{lemma}
\noindent Proof. Under the assumptions of this lemma, a path of type is illustrated in Fig.~\ref{Fig:CCC_external_tangent}. Then, there exists an shorter path from the initial configuration to the target circle, completing the proof. 
$\square$

{Lemmas \ref{LE:3}--\ref{LE:7} indicate that some types in $\mathcal{F}$ cannot be optimal for specific problem parameters. Thus, the results in Lemmas \ref{LE:3}--\ref{LE:7} allow computing the solution without checking all the 12 types in $\mathcal{F}$. 
}




%

\section{Analytical solution for each type in $\mathcal{F}$}\label{SE:Analytical}

Once the final heading angle $\theta_f$ is known, {there are two  positions on the target circle.} The shortest Dubins path from an initial configuration to each of the two positions can be computed in a constant time according to the geometric method in \cite{Dubins:57}. To this end, it amounts to finding the final heading angle in order to address the OCP.

 In this section, based on the geometric properties established in Lemmas \ref{LE:S-pass-center}--\ref{LE:CCC-Straight},  we shall present the relation between the final heading angle $\theta_f$ and the known variables ($x_0$, $y_0$, $\theta_0$, $\rho$, and $r$), so that  $\theta_f$ can be computed analytically. 


\subsection{Analytical solutions for CSC}\label{Subsection:CSC}

For a solution path of $C_1S_2C_3$, let $\alpha_s \in [0,2\pi]$ be the orientation angle of $S_2$, and let $\alpha_c > 0$ be the radian the $C_3$. Then, if $C_3 = R$ (resp. $=L$), we have that the final heading angle is 
 \begin{align}
 \theta_f = \alpha_s - \alpha_c \ \text{(resp.}\ \alpha_s + \alpha_c\text{)}
 \label{EQ:thf_csc}
 \end{align}
where the expression of $\alpha_c $ is given in Lemma \ref{LE:alpha} and $\alpha_s$ can be obtained by simple geometric analysis \cite{Manyam:19}.  It  should be noted from Eq.~(\ref{EQ:thf_csc}) that the length of LSL (resp. RSR) is equal to that of LSR (resp. RSL). Hence, if the solution path of the OCP is of  CSC, then the solution path is not unique.

\subsection{Analytical solution for CCC}

For a solution path of $C_1C_2C_3$, we denote by $[x_c,y_c]$ the center of $C_1$ hereafter. Then, if $C_1 = R$, we have 
$$x_c = x_0 + \rho \cos(\theta_0 -\pi/2)\ \text{and}\ y_c = y_0 + \rho \sin(\theta_0 - \pi/2),$$
and if $C_1 = L$, we have
$$x_c = x_0 + \rho \cos(\theta_0  + \pi/2)\ \text{and}\ y_c = y_0 + \rho \sin(\theta_0 + \pi/2).$$
By the following lemma, we shall show that  the final heading angle $\theta_f$ can be found by solving a quadratic polynomial if the solution path of OCP is of type CCC.
\begin{lemma}\label{LE:CCC-Analytic}
Given any $r>0$, $\rho>0$, and $(x_0,y_0,\theta_0)$, assume  the solution path of the OCP is of type $C_1C_2C_3$. Then, if $C_1C_2C_3 = RLR$ (resp. $= LRL$), we have that $\tan(\frac{\theta_f}{2} -\frac{\pi}{4})$ (resp. $\tan(\frac{\theta_f}{2} + \frac{\pi}{4})$) is a zero of the   quartic polynomial:
\begin{align}
A_1 x^4 + A_2 x^3 + A_3 x^2 + A_4 x + A_5 = 0
\label{EQ:polynomial}
\end{align}
where $A_i$'s are given in Appendix.
\end{lemma} 
Since the proof of this lemma involves some basic mathematical operations, it is delayed to appendix.

{
In summary, since the roots of a quartic polynomial can be  found  either by radicals or by a standard polynomial solver, it follows from Lemma \ref{LE:CCC-Analytic}  that the solution of OCP can be obtained by finding zeros of  a quadratic polynomial   if its type is of CCC. If the solution path is of type CSC, it can be computed  by solving Eq.~(\ref{EQ:thf_csc}). 
}

\section{Numerical simulations}\label{SE:Numerical}

In this section, we present some numerical simulations to illustrate the developments of this paper.

\subsection{Computational cost}

A straightforward way to compute the solution path of the OCP is to  uniformly discretize  the angular position of the target circle.
Given any final angle $\theta_f$, the final position on the target circle is expressed as 
$$\boldsymbol{T}(\delta)\triangleq[r \cos(\theta_f + \delta \pi/2),r \sin (\theta_f + \delta \pi/2)]^T$$ 
where $\delta = 1$ (resp. $\delta = -1$) if the rotational direction of $\theta_f$ is clockwise (resp. counter-clockwise) with respect to center of the target circle. Denote by $D(\theta_f,\delta)$ the shortest Dubins path from  $(x_0,y_0,\theta_0)$ to $\boldsymbol{T}(\delta)$. Then, if the discretization level is denoted by $l>0$,  the Discretization-Based Method (DBM)  is to select an angle $\theta$  in $\{\theta = 2\pi\times i/l: i = 0,1,\ldots,l \}$ so that $D(\theta,\delta)$ is the smallest, i.e.,
\[D(\theta,\delta) = \underset{i=0,1,\ldots,l}{\min}{D(2\pi\times i/l,\delta)},\ \ \ \delta = \pm 1.\]
For notational simplicity, we denote hereafter by DBM($l$) the DBM with a discretization level of $l\in \mathbb{N}$.

Let the parameters $(x_0,y_0,\theta_0)$, $\rho>0$, and $r>0$ be generated randomly by uniform distribution. Both   the analytical solutions in Section \ref{SE:Analytical} and the DBM(360) are tested on 10000 randomly generated OCPs.   The computational time of the analytic method is tested by MATLAB on a desktop with Intel(R) Core(TM)i7-8550U CPU @1.80 GHz, in comparison with the DBM(360). Table \ref{Tab:Compare} presents the computation time for different $\rho$, where $d_m\geq 0$ denotes the distance from initial point $(x_0,y_0)$ to the center of the target circle. 
\begin{table*}
\caption{Time consumed by the analytic method  and the DBM(360)}\label{Tab:Compare}
\centering
\begin{tabular}{c|ccccc}
\hline
\diagbox{Method}{Consumed Time (s)}{$d_m$} & $\geq 4\rho$ & $= 3\rho$ & $=2\rho$ & $=\rho$  & $<\rho$\\
\hline
Analytic &  1.15$\times 10^{-4}$ & 1.34$\times 10^{-4}$ & 1.39$\times10^{-4}$ & 1.31$\times 10^{-4}$ &1.48$\times 10^{-4}$\\
\hline
 DBM(360) & 0.2386 & 0.2276 & 0.2808  & 0.2836 & 0.2920\\
\hline
\end{tabular}
\end{table*}
We can see from Table \ref{Tab:Compare} that the improvement factors of the analytic method compared with the DBM(360) are greater than around 2000. 

Note that the DBM(360) can only generate an approximate solution  for the OCP. If a more accurate solution is required, a higher level of discretization is needed, which however results in a higher computational cost. As the analytic solution to each type in $\mathcal{F}$ has been devised, the accurate solution of the OCP can be obtained in a constant time by checking some analytic  functions.

\subsection{Specific examples}

In this subsection, we present some  examples to illustrate the developments of the paper.

\subsubsection{Case A}
\begin{figure}[!htp]
\centering
\begin{subfigure}[t]{4cm}
\centering
\includegraphics[width = 4cm]{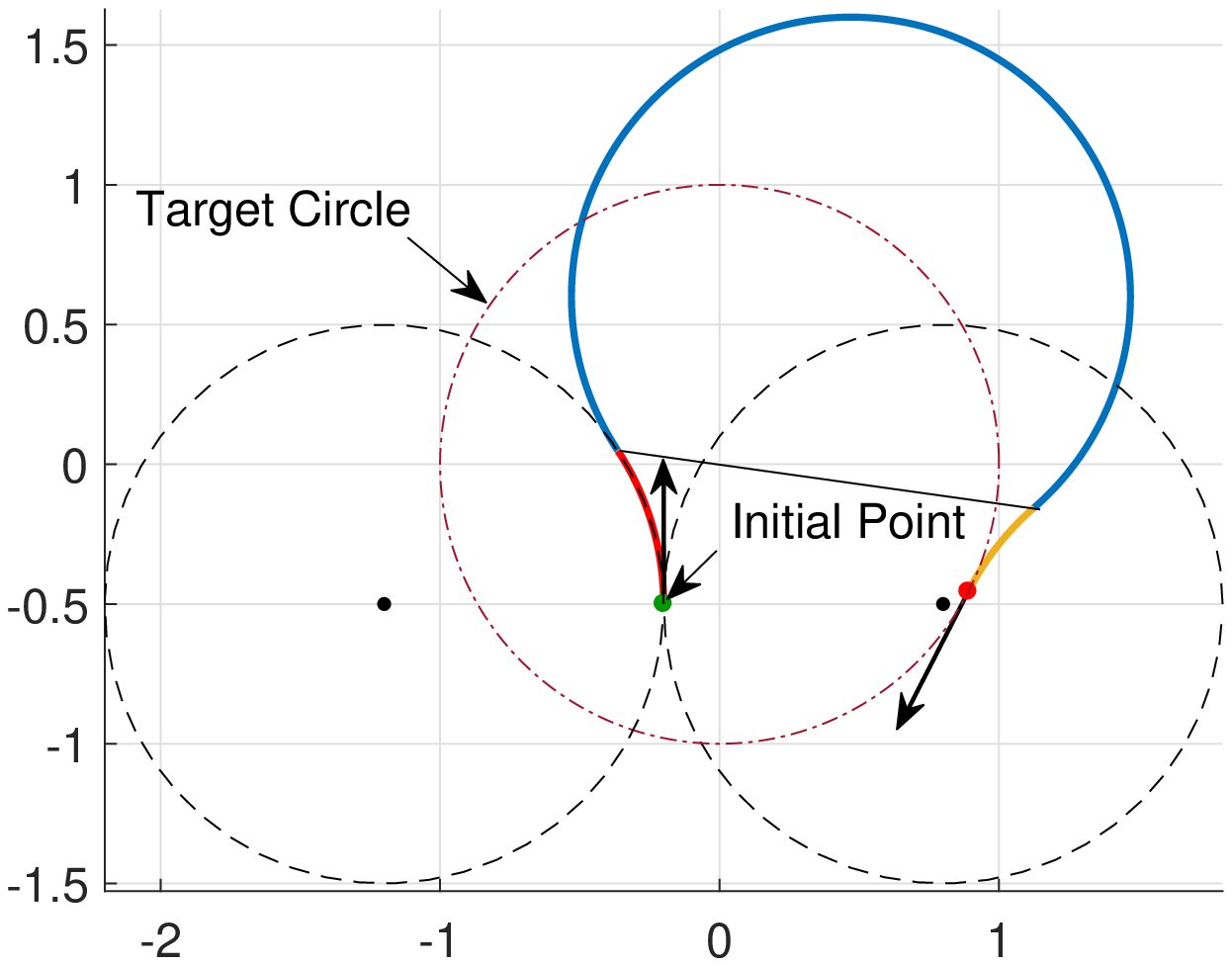}
\caption{Clockwise}
\label{Fig:caseB1}
\end{subfigure}
\begin{subfigure}[t]{4cm}
\centering
\includegraphics[width = 4cm]{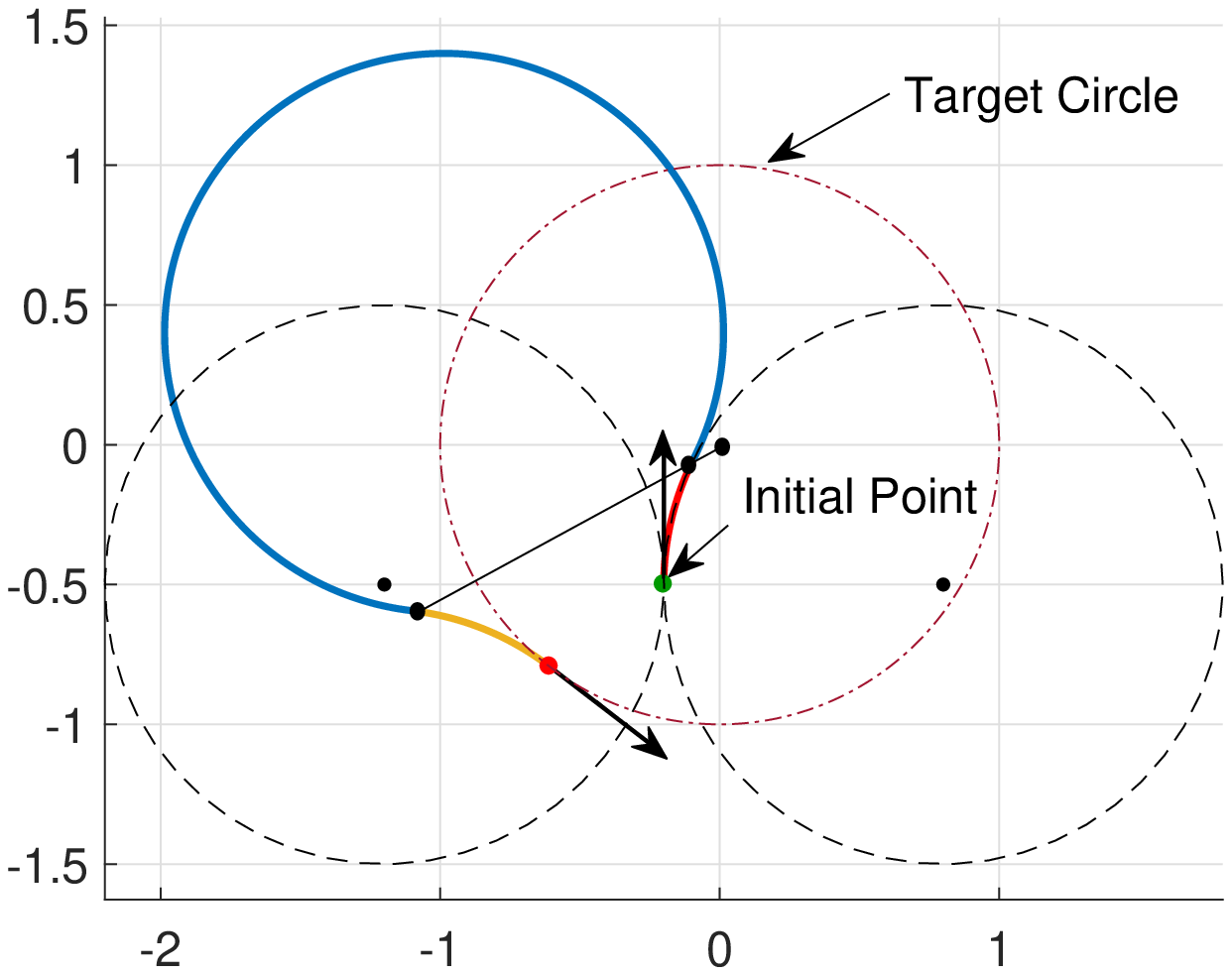}
\caption{Counter-clockwise}
\label{Fig:caseB2}
\end{subfigure}
\caption{The solution paths of the OCP for different rotational directions for case A.}
\label{Fig:CaseB}
\end{figure}

Set $\rho = 1$, $r = 1$, and $(x_0,y_0,\theta_0) = (-0.2,-0.5,\pi/2)$.  The analytical results in Section \ref{SE:Analytical} are  applied to computing the shortest path from $(x_0,y_0,\theta_0)$ to the target circle.  Fig.~\ref{Fig:CaseB} shows the solution paths with two different final rotational directions.
It is apparent to see from Fig.~\ref{Fig:CaseB} that the two concatenating points and  center of the target circle are collinear, coinciding with Lemma \ref{LE:CCC-Straight}.

\subsubsection{Case B}

\begin{figure}[!htp]
\centering
\includegraphics[width = 4cm]{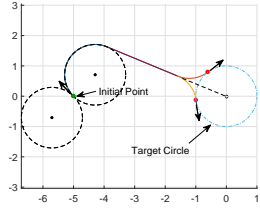}
\caption{The shortest paths for case B.}
\label{Fig:caseB-B}
\end{figure}

The initial condition is set as $(x_0,y_0,\theta_0 ) = (-5,0,3\pi/2)$. Set $\rho =1 $ and $ r = 1$.
 The shortest paths are computed by directly applying  the analytical results in Section \ref{SE:Analytical} and presented in Fig.~\ref{Fig:caseB-B}.
We can see from Fig.~\ref{Fig:caseB-B} that the straight line segment and  center of the target circle are collinear, as predicted by Lemma \ref{LE:S-pass-center}. Also notice that the length of the path of RSR is the same as that of the path of RSL, as analyzed in Subsection \ref{Subsection:CSC}.

\subsubsection{Case C}

\begin{figure}[!htp]
\centering
\includegraphics[width = 4cm]{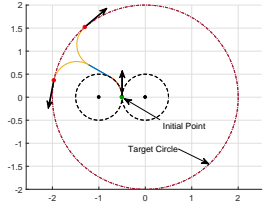}
\caption{The shortest paths for case C.}
\label{Fig:caseC}
\end{figure}

Set $r = 2$, $\rho = 0.5$, and $(x_0,y_0,\theta_0) = (-0.5,0,\pi/2)$. Employing the analytical results in Section \ref{SE:Analytical} once again, the shortest paths are computed and presented in Fig.~\ref{Fig:caseC}. 
We can see that the final circular arc is internally tangent to the target circle. In this case, we still have that the straight line segment passes is collinear with the center of the target circle, as predicted by Lemma \ref{LE:S-pass-center}, and that the length of the path of LSL is the same as that of the path of LSR.

\subsubsection{Case D}

\begin{figure}[!htp]
\centering
\begin{subfigure}[t]{4cm}
\centering
\includegraphics[width = 4cm]{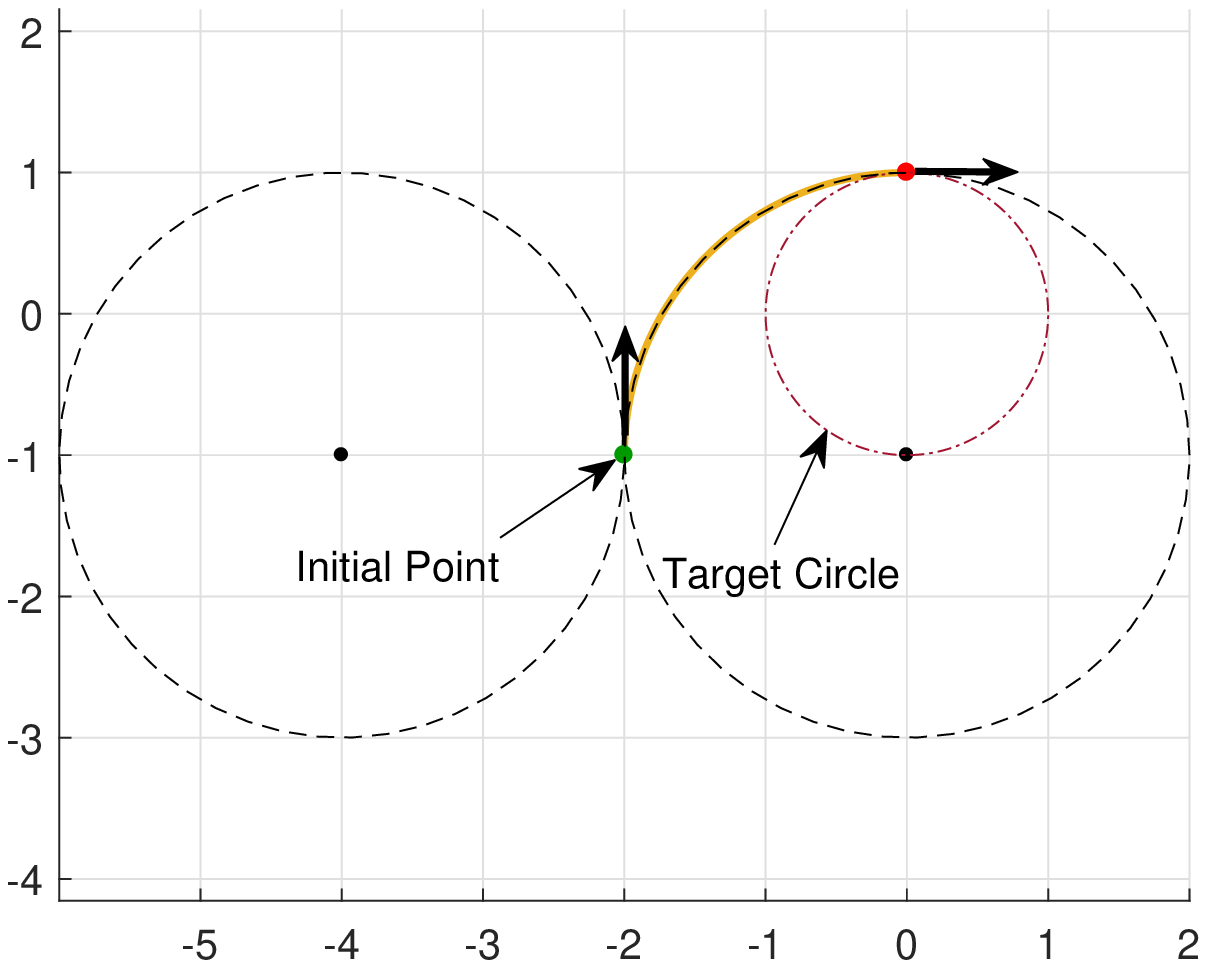}
\caption{Analytic method}
\label{Fig:caseD1}
\end{subfigure}
\begin{subfigure}[t]{4cm}
\centering
\includegraphics[width =4cm]{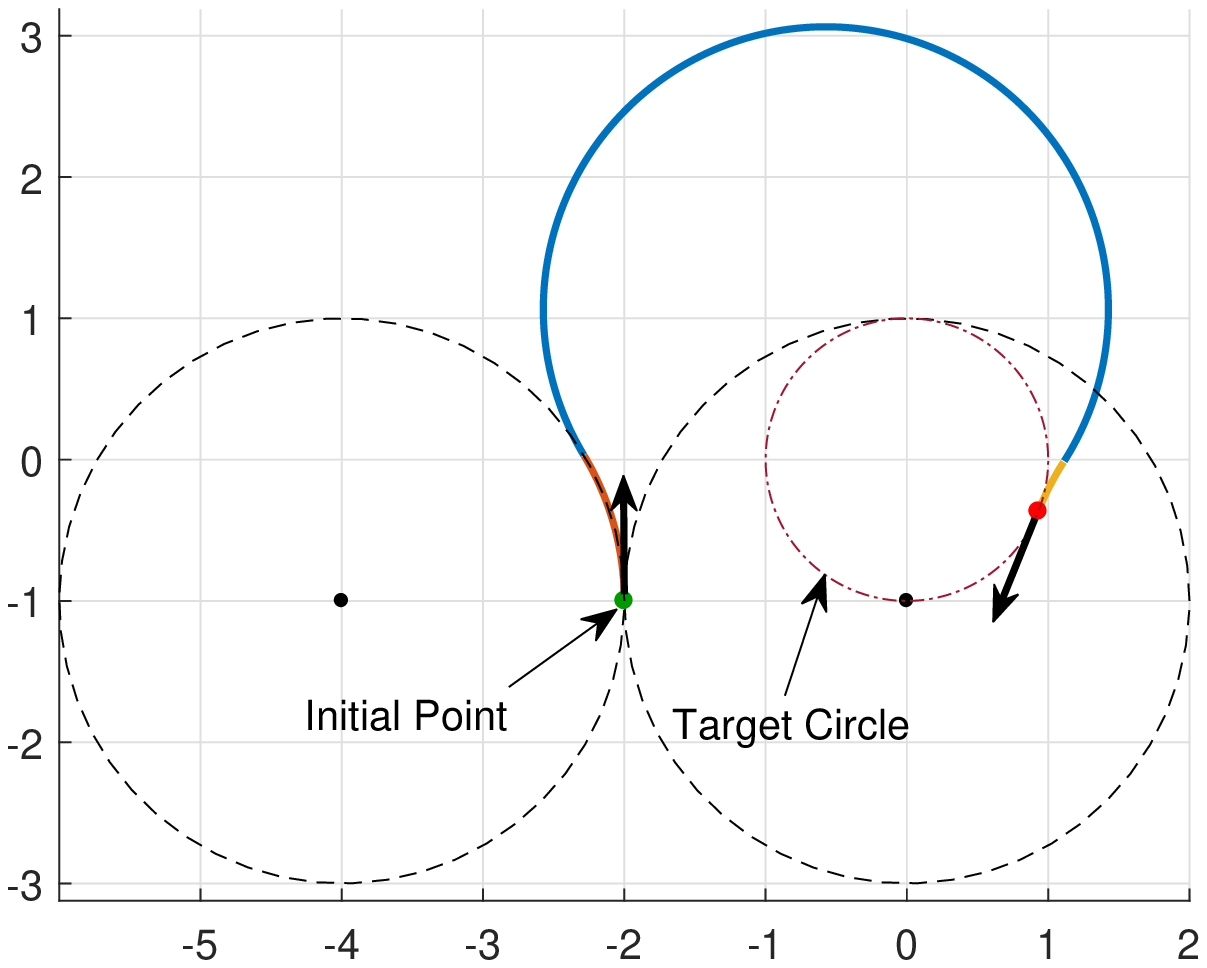}
\caption{DBM(360)}
\label{Fig:caseD2}
\end{subfigure}
\caption{The solution paths of the OCP for case D.}
\label{Fig:caseD}
\end{figure}

Let $r = 1$ and $\rho = 2$, and set the initial condition as 
\begin{align}
\left[
\begin{array}{c}
x_0\\
y_0\\
\theta_0
\end{array}
\right]
= 
\left[
\begin{array}{l}
r \cos (269.5*\pi/180) + \rho \cos (179.5*\pi/180) \\
r \sin(269.5*\pi/180) + \rho \sin (179.5*\pi/180)\\
89.5*\pi/180
\end{array}
\right]\nonumber
\end{align}
This initial condition is tailored so that the solution path of the OCP is a single circular arc.  Using the analytical solutions in Section \ref{SE:Analytical}, the solution is computed and presented in Fig.~\ref{Fig:caseD1}.
However, it should be noted that the DBM(360) cannot find the solution path. In fact, the path computed by the DBM(360) is quite different from the analytical solution path, as shown in Fig.~\ref{Fig:caseD2}. To compute accurate solution for case D by DBM($l$), the discretization level $l$ should be large enough, which however will result in higher computational time (note that the computational time of DBM($l$) is in direct proportion to the discretization level $l$). 

Taking into account all the numerical examples presented in this section, it is concluded that the analytic method not only can compute the solution path of the OCP in a constant time but also can generate more accurate solutions, in comparison with the DBM(360).

\section{Conclusions}

The shortest Dubins paths  from a fixed configuration to a target circle with the terminal heading (or velocity) tangential to the circle was studied by applying PMP.  
Through synthesizing  the necessary conditions for optimality, some geometric properties for the shortest path were presented. To be  specific, once the shortest path is of type CSC,   the straight line segment S is collinear with the center of the target circle; if the shortest path is of type CCC, then the two concatenating points between the circular arcs and the center of the target circle are collinear. These geometric properties  ruled out the substring CS so that the shortest path must lie in a sufficient family of 12 types. The geometric properties allowed to devise an analytical solution for each of all the 12 types. In addition, some relationships between problem parameters and geometric properties were revealed, allowing finding the solution path without checking all the 12 types.  Comparing with the straightforward discretization-based method,  the developments in the paper not only enabled reducing computation time (cf. Table \ref{Tab:Compare}) but also allowed generating more accurate solutions.

\section*{Acknowledgments}
This research was   supported by the National Natural Science Foundation of China (No. 61903331).

\section*{Appendix: Proof of Lemma \ref{LE:CCC-Analytic}}

\begin{figure}[!htp]
\centering
\begin{subfigure}[t]{4cm}
\centering
\includegraphics[width = 4cm]{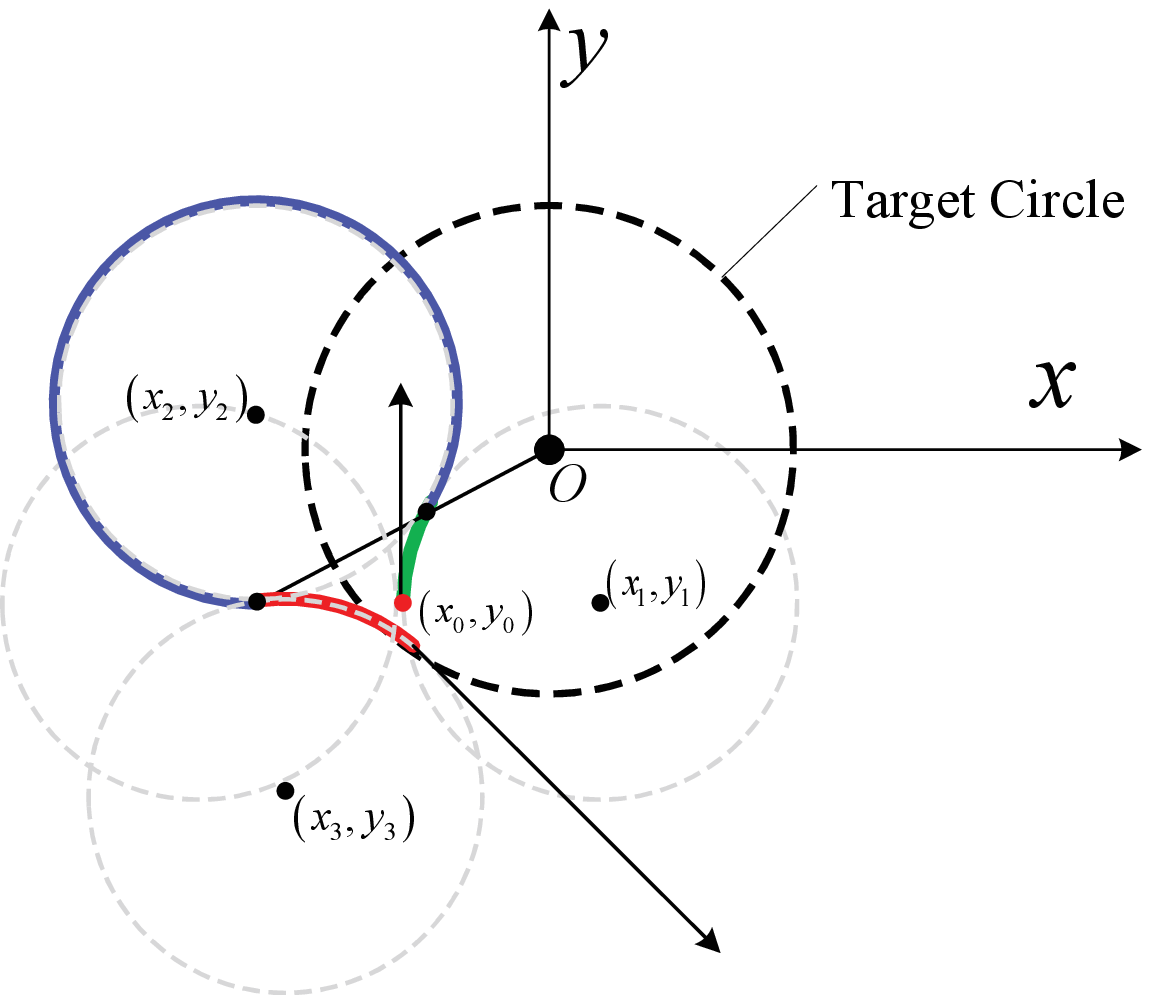}
\caption{RLR}
\label{Fig:RLR}
\end{subfigure}
\begin{subfigure}[t]{4cm}
\centering
\includegraphics[width = 4cm]{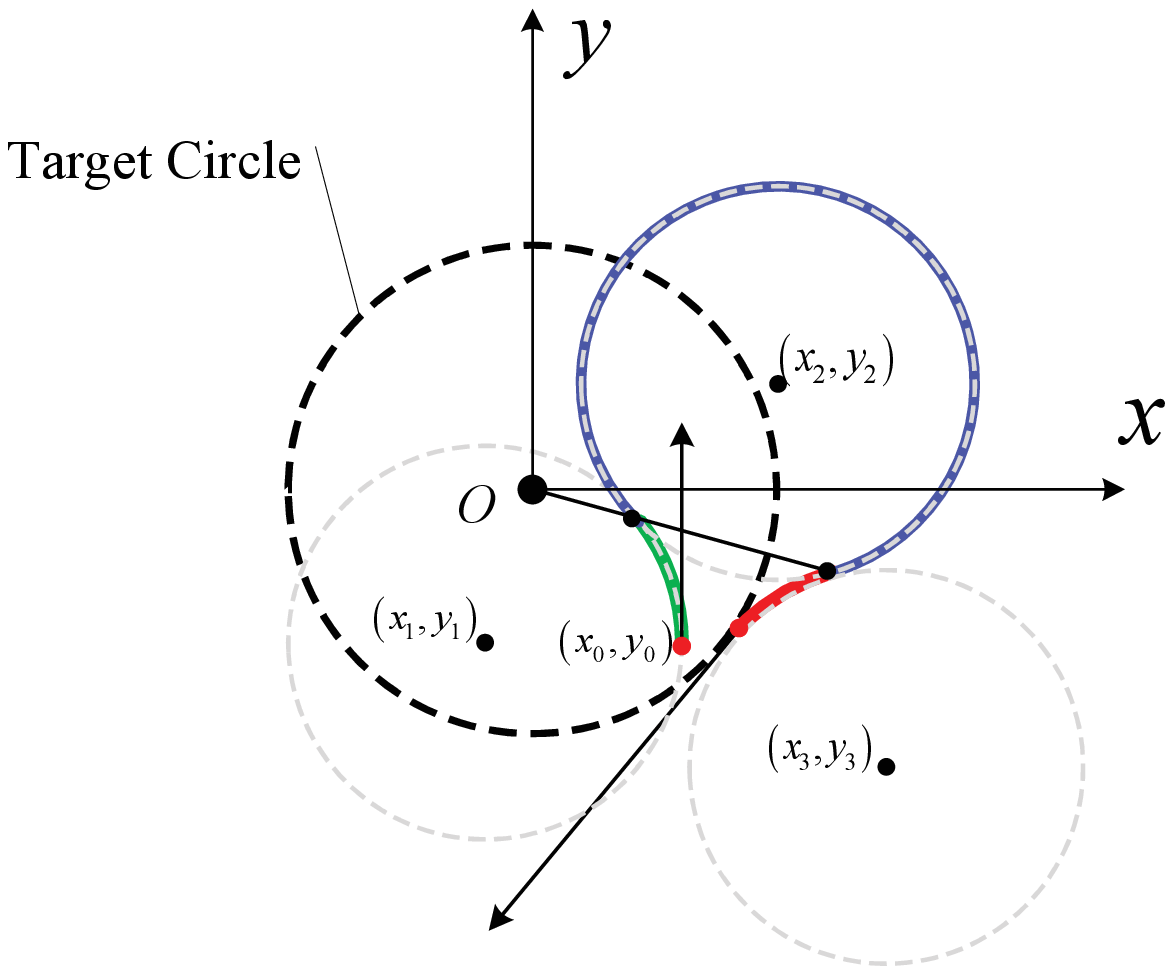}
\caption{LRL}
\label{Fig:LRL}
\end{subfigure}
\caption{The geometry for the RLR and LRL paths.}
\label{Fig:CCC}
\end{figure}

First consider $C_1C_2C_3= RLR$, as shown by Fig.~\ref{Fig:RLR} where the circular arc $C_3$ is externally tangent to the target circle.  
Denote by $\boldsymbol{c}_3^r  = [x_3,y_3]^T$ the center of  $C_3$. According to Lemma \ref{LE:CCC-ruleout-internal}, if   $C_3$ is internally tangent to the target circle, we have $\rho < r$. 
Hence, the center $\boldsymbol{c}_3^r $ lies on the circle with radius $  r + \delta \rho$ and centered at the origin where $\delta = 1$ (resp. $\delta = -1$) if the target circle and $C_3$ are externally (resp. internally) tangent.  Note that the vector from the origin  to $\boldsymbol{c}_3^r$ is perpendicular with the final heading, indicating that the vector $[\cos(\theta_f -\pi/2),\sin(\theta_f - \pi/2)]^T$ is align to the vector $[x_3,y_3]^T$.  Therefore, $\c_3^r$ explicitly writes
\begin{align}
\left[\begin{array}{c}
x_3\\
y_3
\end{array}
\right]  = (r+ \delta \rho)
\left[\begin{array}{c} \cos (\theta_f - \pi/2)\\
  \sin (\theta_f - \pi/2)\end{array}
\right].
\label{EQ:lemma-ccc-3}
\end{align} 
Set $\boldsymbol{c}_1^r: = [x_c,y_c]^T$.
Denote by $\boldsymbol{c}_2^l = [x_2,y_2]^T$ the center of the second circular arc $C_2$, by $A = [x_A ,y_A]^T$ and $B = [x_B,y_B]^T$ the concatenating points from $C_1$ to $C_2$ and from $C_2$ to $C_3$, respectively. Furthermore, denote by $\alpha_A$ and $\alpha_B$ the heading angle at $A$ and $B$, respectively.  Then, because the vector from $\c_1^r$ to $\c_2^l$ is perpendicular with the heading at $A$ and the vector from $\c_3^r$ to $\c_2^l$ is perpendicular with the heading at $B$,  we have
\[\boldsymbol{c}_2^l = \boldsymbol{c}_1^r + 2\rho[\cos (\alpha_A + \pi/2),\sin(\alpha_A + \pi/2)]^T\] and 
\[\boldsymbol{c}_2^l  = \boldsymbol{c}_3^r + 2\rho[\cos (\alpha_B + \pi/2),\sin(\alpha_B + \pi/2)]^T \]
Combining these two equations indicates
\begin{align}
x_3 + 2 \rho \cos (\alpha_B + \pi/2) = x_c + 2\rho \cos (\alpha_A + \pi/2)\label{EQ:analytical-00}\\
y_3  + 2\rho  \sin(\alpha_B + \pi/2) = y_c + 2\rho \sin(\alpha_A + \pi/2)\label{EQ:analytical-01}
\end{align}
These two equations can be simplified to
\begin{align}
0 \ & = (x_3 - x_c)^2 + 4(x_3 - x_c)\rho \cos (\alpha_B + \pi/2)\nonumber\\
\ &  + (y_3 - y_c)^2 + 4 (y_3 - y_c) \rho \sin (\alpha_B+\pi/2)
\label{EQ:analytical_10}
\end{align}
According to the geometry in Fig.~\ref{Fig:RLR}, we also have
\begin{align}
\left[
\begin{array}{c}x_A\\
y_A
\end{array}\right] = \left[
\begin{array}{c} x_c
\\
y_c
\end{array}\right]
+ \rho
\left[\begin{array}{c} \cos (\alpha_A + \pi/2) \\
 \sin(\alpha_A + \pi/2)
\end{array}\right] 
\label{EQ:xA}\\
\left[
\begin{array}{c}x_B\\
y_B
\end{array}\right] = \left[
\begin{array}{c} x_3 
\\
y_3\end{array}\right]
+ \rho
\left[
\begin{array}{c} \cos (\alpha_B + \pi/2) \\
\sin(\alpha_B + \pi/2)
\end{array}\right]
\label{EQ:xB}
\end{align}
Since $A$, $B$, and the origin $O$ lie on a straight line according to Lemma \ref{LE:CCC-Straight}, it follows 
\begin{align}
y_A / x_A = y_B / x_B.
\label{EQ:AB}
\end{align}
Substituting Eq.~(\ref{EQ:xA}) and Eq.~(\ref{EQ:xB}) into Eq.~(\ref{EQ:AB}) leads to
\begin{align}
\frac{y_c + \rho \sin(\alpha_A + \pi/2)}{x_c + \rho \cos (\alpha_A + \pi/2)  } = \frac{y_3 + \rho \sin(\alpha_B + \pi/2)}{ x_3 + \rho \cos (\alpha_B + \pi/2)}.
\label{EQ:CCC-analytical-20}
\end{align}
Substituting Eq.~(\ref{EQ:analytical-00}) and Eq.~(\ref{EQ:analytical-01}) into Eq.~(\ref{EQ:CCC-analytical-20}), we have
\begin{align}
\frac{y_c + y_3 + 2\rho \sin (\alpha_B + \frac{\pi}{2})}{x_c + x_3 + 2\rho \cos (\alpha_B +  \frac{\pi}{2})}= \frac{y_3 + \rho \sin (\alpha_B +  \frac{\pi}{2})}{x_3 + \rho \cos (\alpha_B +  \frac{\pi}{2})}
\end{align}
which can be rearranged as
\begin{align}
0 \ & =
y_c x_3 - x_c y_3 + (x_3 - x_c)\rho \sin (\alpha_B + \pi/2)\nonumber\\
\ & +(y_c - y_3)\rho \cos (\alpha_B + \pi/2)
\label{EQ:analytical_11}
\end{align}
Combining Eq.~(\ref{EQ:analytical_10}) with Eq.~(\ref{EQ:analytical_11}) leads to
\begin{align}
0 \ & =(x_c y_3 - y_c x_3)^2 + [(x_3 - x_c)^2 + (y_3 - y_c)^2]^2/16 \nonumber\\
& -  \rho^2[(x_3 - x_c)^2 + (y_3 - y_c)^2] .
\label{EQ:temp}
\end{align}
Substituting Eq.~(\ref{EQ:lemma-ccc-3}) into this equation yields
\begin{align}
0&\ = (r+\delta \rho)^2 [x_c \sin(\theta_f -\pi/2) - y_c \cos(\theta_f - \pi/2)]^2\nonumber\\
&\ + \{- 2  (r + \delta \rho)[ x_c \cos(\theta_f - \pi/2)+ y_c \sin (\theta_f -\pi/2)]\nonumber\\
&\  + (r+\delta \rho)^2 +  x_c^2 + y_c^2 \}^2/16 - \rho^2 \{ x_c^2  + y_c^2 + (r+\delta \rho)^2  \nonumber\\
&\ - 2  (r + \delta \rho)[ x_c \cos(\theta_f - \pi/2)+ y_c \sin (\theta_f -\pi/2)]  \}
\label{EQ:lemma-RLR-poly0}
\end{align}
By taking into account the half-angle formulas
\begin{align}
\sin \theta = \frac{2 \tan(\frac{\theta}{2})}{1 + \tan^2 (\frac{\theta}{2})} \ \text{and}\ 
\cos \theta = \frac{1 -  \tan^2(\frac{\theta}{2})}{1 + \tan^2 (\frac{\theta}{2})}\nonumber
\end{align}
 Eq.~(\ref{EQ:lemma-RLR-poly0}) can be written as the quadratic polynomial  in Eq.~(\ref{EQ:polynomial}) whose coefficients are expressed as
 \begin{align}
A_1 & =  2 x_c^2 \left(\rho ^2+9 r^2+18\delta \rho  r-2 y_c (\delta \rho +r)+y_c^2\right ) + \nonumber\\
& \left(\delta \rho +r+y_c\right){}^2 \left[(4+\delta) \rho +r+y_c\right]\times\nonumber\\
&  \left[(-4 + \delta) \rho +r+y_c\right]+x_c^4\nonumber\\
A_2 & = -[-7 \rho ^2+r^2+2 \delta  \rho  r-6 y_c (\delta  \rho +r)+x_c^2+y_c^2]\nonumber\\
& \times 8 x_c (\delta  \rho +r)  \nonumber\\
A_3 & =2 [-2 x_c^2 \left(11 \rho ^2+3 r^2+6 \delta  \rho  r\right) +x_c^4+y_c^4\nonumber\\
&  +2 y_c^2 \left(7 \rho ^2+15 r^2+30 \delta  \rho  r+x_c^2\right)+\nonumber\\
& ((\delta -4) \rho +r) (\delta  \rho +r)^2 ((\delta +4) \rho +r) ] \nonumber\\
A_4 & =  \left(-7 \rho ^2+r^2+2 \delta  \rho  r-6 y_c (\delta \rho +r)+x_c^2+y_c^2\right)\times\nonumber\\
& [-8 x_c (\delta \rho +r)]\nonumber\\
A_5 & = 2 x_c^2 \left(\rho ^2+9 r^2+18 \delta \rho  r+2 y_c (\delta \rho +r)+y_c^2\right)+\nonumber\\
& \left(\delta \rho +r-y_c\right){}^2 \left((4+ \delta) \rho +r-y_c\right)\times \nonumber\\
& \left((-4 + \delta)\rho +r-y_c\right)+x_c^4\nonumber
\end{align}
with $\delta = 1$ (resp. $\delta = -1$) if  $C_3$ is externally (resp. internally) tangent to the target circle.  For $C_1C_2C_3= LRL$, it can be proved in the same way, so the proof for the type of LRL is omitted here.



\bibliographystyle{plain}        
\bibliography{autosam}

\end{document}